\documentclass[11pt,leqno]{article}
\usepackage{amsthm,amsfonts,amssymb,amsmath,oldgerm}
\usepackage{epsfig}
\numberwithin{equation}{section} 

\renewcommand\d{\partial}

\def\l{\lambda}
\def\eps{\varepsilon }
\def\e{\varepsilon} 

\newcommand{\rvec}{{\cal R}}

\renewcommand\d{\partial}

\newcommand\R{\mathbb R}
\newcommand\C{\mathbb C}

\def\eps{\varepsilon}
\def\e{\varepsilon}

\def\l{\lambda}

\newcommand\br{\begin{remark}}
\newcommand\er{\end{remark}}
\newcommand\bp{\begin{pmatrix}}
\newcommand\ep{\end{pmatrix}}
\newcommand\be{\begin{equation}}
\newcommand\ee{\end{equation}}
\newcommand\ba{\begin{equation}\begin{aligned}}
\newcommand\ea{\end{aligned}\end{equation}}

\newcommand{\bap}{\begin{app}}
\newcommand{\eap}{\end{app}}
\newcommand{\begs}{\begin{exams}}
\newcommand{\eegs}{\end{exams}}
\newcommand{\beg}{\begin{example}}
\newcommand{\eeg}{\end{exaplem}}
\newcommand{\bpr}{\begin{proposition}}
\newcommand{\epr}{\end{proposition}}
\newcommand{\bt}{\begin{theorem}}
\newcommand{\et}{\end{theorem}}
\newcommand{\bc}{\begin{corollary}}
\newcommand{\ec}{\end{corollary}}
\newcommand{\bl}{\begin{lemma}}
\newcommand{\el}{\end{lemma}}
\newcommand{\bd}{\begin{definition}}
\newcommand{\ed}{\end{definition}}
\newcommand{\brs}{\begin{remarks}}
\newcommand{\ers}{\end{remarks}}

\newtheorem{theo}{Theorem}[section]
\newtheorem{prop}[theo]{Proposition}

\newtheorem{exam}[theo]{Example}

\newtheorem{exams}[theo]{Examples}

\numberwithin{equation}{section}

\newcommand{\CalF}{\mathcal{F}}
\newcommand{\CalC}{\mathcal{C}}

\newcommand{\RR}{{\mathbb R}}

\newcommand{\CC}{{\mathbb C}}

\newcommand{\Span}{{\rm Span }}

\newtheorem{theorem}{Theorem}[section]
\newtheorem{proposition}[theorem]{Proposition}
\newtheorem{corollary}[theorem]{Corollary}
\newtheorem{lemma}[theorem]{Lemma}
\newtheorem{definition}[theorem]{Definition}

\newtheorem{example}[theorem]{Example}
\newtheorem{remark}[theorem]{Remark}

\newcommand\cN{{\cal  N}}
\newcommand\cE{{\cal  E}}

\newcommand\cM{{\mathcal M}}

\newcommand\tG{{\tilde G}}

\title{
The refined inviscid stability condition and\\ 
cellular instability of viscous shock waves
}

\author{\sc \small 
Kevin Zumbrun\thanks{Indiana University, Bloomington, IN 47405;
kzumbrun@indiana.edu:
Research of K.Z. was partially supported
under NSF grants no. DMS-0300487 and DMS-0801745.
 }}
\begin{document}

\maketitle


\begin{abstract}
Combining work of 
Serre and Zumbrun, Benzoni-Gavage, Serre, and Zumbrun, and Texier and Zumbrun,
we propose as a mechanism for the onset of cellular instability
of viscous shock and detonation
waves in a finite-cross-section duct the violation of
the refined planar stability condition of Zumbrun--Serre,
a viscous correction of the inviscid planar stability condition of Majda.
More precisely, we show for a model problem involving flow 
in a rectangular duct with artificial periodic boundary conditions
that transition to multidimensional instability 
through violation of the refined stability condition of
planar viscous shock waves on the whole space generically implies for a duct
of sufficiently large cross-section a cascade of Hopf bifurcations
involving more and more complicated cellular instabilities.
The refined condition is numerically calculable as described
in Benzoni-Gavage--Serre-Zumbrun.
\end{abstract}

\bigbreak
\section{Introduction}

It is well known both experimentally and numerically
\cite{BE, MT, BMR, FW, MT, ALT, AT, F1, F2,KS}
that shock and detonation waves propagating in a finite cross-section
duct can exhibit time-oscillatory or ``cellular'' instabilities,
in which the initially nearly planar shock takes on nontrivial 
transverse geometry.
Majda et al \cite{MR1,MR2,AM} have
studied the onset of such instabilities by weakly nonlinear optics
expansion of the associated planar inviscid shock in the whole space.
More recently, Kasimov--Stewart \cite{KS} and Texier--Zumbrun \cite{TZ2,TZ3,TZ3}
have studied these instabilities as Hopf bifurcations of flow
in a finite-cross-section duct, associated
with passage across the imaginary axis of eigenvalues of the linearized
operator about the wave.

In this paper, combining the analyses of \cite{BSZ,TZ4,Z4},
we make an explicit connection between stability
of planar shocks on the whole space, and Hopf bifurcation in
a finite cross-section duct, by a mechanism different from that
investigated by Majda et al.
Specifically, we point out that violation of
 the {\it refined stability condition} of \cite{ZS,Z1,BSZ},
a viscous correction of the inviscid planar stability condition of Majda
\cite{M1}--\cite{M4}, is generically associated with Hopf bifurcation
in a finite cross-section duct corresponding to the observed cellular
instability, for cross-section $M$ sufficiently large.
Indeed, we show more, that this is associated with a cascade of
bifurcations to higher and higher wave numbers and more and more
complicated solutions,
with features on finer and finer length/time scales.

\subsection{Equations and assumptions}
Consider a planar viscous shock solution
\be\label{planar}
u(x,t)=\bar u(x_1-st)
\ee
of a two-dimensional system of viscous conservation laws
\be\label{multiNS}
u_t + \sum f^j(u)_{x_j}= \Delta_x u, 
\qquad u\in \RR^n, \, x\in \RR^2, \, t\in \RR^+
\ee
on the whole space. 
This may be viewed alternatively as a planar traveling-wave solution 
of \eqref{multiNS} on an infinite channel
$$
\CalC:= \{x:\, (x_1, x_2)\in \RR^1\times [-M,M]\}
$$
under periodic boundary conditions
\be\label{BC}
u(x_1,M )=u(x_1,-M).
\ee

We take this as a simplified mathematical model for compressible flow in
a duct, in which we have neglected boundary-layer phenomena along the
wall $\partial \Omega$ in order to isolate the oscillatory phenomena
of our main interest.

Following \cite{TZ2},
consider a one-parameter family of standing planar viscous shock 
solutions $\bar u^\eps(x_1)$ of a smoothly-varying family of conservation laws 
\begin{equation}
\label{multisysepseqn}
u_t =\CalF(\e, u):= \Delta_x u- \sum_{j=1}^2 F^j(\e, u)_{x_j},
\qquad u\in \RR^n
\end{equation}
in a fixed channel $\CalC$, with 
periodic boundary conditions (typically, shifts 
$\sum F^j(\eps, u)_{x_j}:= \sum f^j(u)_{x_j} -s(\eps) u_{x_1}$
of a single equation \eqref{multiNS}
written in coordinates
$x_1\to x_1-s(\e)t$ moving with traveling-wave solutions
of varying speeds $s(\e)$),
with linearized operators
$L(\e) :=\partial \CalF/\partial u|_{u=\bar u^\e}$.
Profiles $\bar u^\eps$ satisfy the standing-wave ODE
\be\label{multiode}
u'=F^1(\e,u)- F^1(\e,u_-).
\ee

Let 
\begin{equation} \label{A'}
 A^1_\pm(\e) := \lim_{z \to \pm \infty} F^1_u(\e,\bar u^\e).
 \end{equation}
Following \cite{Z1,TZ2,Z4}, we make the assumptions:

\quad (H0) \quad  $F^j\in C^{k}$, $k\ge 2$.
\medbreak
\quad (H1)  \quad 
$\sigma (A^1_\pm(\eps))$ real, distinct, and nonzero,
and $\sigma (\sum \xi_j A^j_\pm(\eps))$ real and semisimple for $\xi\in \RR^d$.
\medbreak

\noindent For most of our results, we require also:

\quad (H2)  \quad  Considered as connecting orbits of \eqref{multiode}, 
$\bar u^\eps$ are transverse and unique up to translation,
with dimensions of the stable subpace $S(A^1_+)$ 
and the unstable subspace $U(A^1_-)$ 
summing for each $\eps$ to $n+1$.
\medbreak

\quad (H3)  \quad  
$\det(r_1^-, \dots, r_c^-, r_{c+1}^+, \dots, r_n^+, u_+-u_-)\ne 0$,
where $r_1^-,\dots, r_c^-$ are eigenvectors of $A^1_-$ associated
with negative eigenvalues and
$r_{c+1}^-,\dots, r_n^-$ are eigenvectors of $A^1_+$ associated
with positive eigenvalues.
\medbreak

Hypothesis (H2) asserts in particular that
$\bar u^\eps$ is of standard {\it Lax type}, 
meaning that the axial hyperbolic
convection matrices $A^1_+(\eps)$ and $A^1_-(\eps)$ 
at plus and minus spatial infinity have, respectively, 
$n-c$ positive and $c-1$ negative real eigenvalues for 
$1\le c \le n$, where $c$ is the 
characteristic family associated with the shock: in other words, there
are precisely $n-1$ outgoing hyperbolic characteristics in the far field.
Hypothesis (H3) may be recognized as the {\it Liu--Majda condition}
corresponding to one-dimensional stability of the associated inviscid shock.
In the present, viscous, context, this, together with transverality,
(H2), plays the role of a spectral nondegeneracy condition
corresponding in a generalized sense \cite{ZH,Z1} to simplicity
of the embedded zero eigenvalue associated with eigenfunction
$\partial_{x_1}\bar u$ and translational invariance.

\subsection{Stability conditions}

Our first set of results, generalizing the one-dimensional analysis
of \cite{Z4}, characterize stability/instability of waves
$\bar u^\eps$ in terms of the spectrum of the linearized
operator $L(\eps)$.
Fixing $\eps$, we suppress the parameter $\eps$.
We start with the routine observation
that the semilinear parabolic equation \eqref{multiNS}
has a center-stable manifold about the equilibrium solution $\bar u$.

\bpr\label{t:maincs}
Under assumptions (H0)--(H1), 
there exists in an $H^2$
neighborhood of the set of translates of $\bar u$
a codimension-$p$ translation invariant $C^k$ (with respect to $H^2$) 
center stable manifold $\cM_{cs}$,
tangent at $\bar u$ to the center stable subspace $\Sigma_{cs}$ of $L$, 
that is (locally) invariant under the forward time-evolution of 
\eqref{multiNS}--\eqref{BC} 
and contains all solutions that remain bounded and 
sufficiently close to a translate of $\bar u$ in forward time, where $p$ is the
(necessarily finite) number of unstable, i.e., positive real part,
eigenvalues of $L$.  
\epr

\begin{proof}
By standard considerations \cite{He,TZ2}, 
$L(\eps)$ possesses no essential spectrum and 
at most a finite set of positive real part eigenvalues
on $\Re \lambda>0$.
With this observation, the result
follows word-for-word by the argument of \cite{Z4} in
the one-dimensional case, which depends only on the properties
of $L$ as a sectorial second-order elliptic operator, and on semilinearity
and translation-invariance of the underlying equations \eqref{multiNS}.
\end{proof}

Introduce now the {\it nonbifurcation condition}:

\quad (D1) \quad $L$ has no nonzero imaginary eigenvalues. 

As discussed above, (H2)--(H3) correspond 
to a generalized notion of simplicity of the 
embedded eigenvalue $\lambda=0$ of $L$.
Thus, (D1) together with (H2)--(H3) correspond to the 
assumption that there are no additional (usual or generalized) 
eigenvalues on the imaginary 
axis other than the translational eigenvalue at $\lambda=0$;
that is, the shock is not in transition between different degrees of
stability, but has stability properties that are insensitive to 
small variations in parameters.

\bt\label{t:mainstab}
Under (H0)--(H3) and (D1), $\bar u$ is nonlinearly
orbitally stable as a solution of \eqref{multiNS}--\eqref{BC} 
under sufficiently small perturbations
in $L^1\cap H^2$ lying on the codimension $p$ center stable 
manifold $\cM_{cs}$ of $\bar u$ and its translates, 
where $p$ is the number of unstable eigenvalues of $L$,
in the sense that, for some $\alpha(\cdot)$, all $L^p$,
\ba\label{bounds}
|u(x, t)-\bar u(x-\alpha(t))|_{L^p}&\le
C(1+t)^{-\frac{1}{2}(1-\frac{1}{p})}
|u(x,0)-\bar u(x)|_{L^1\cap H^2},
\\
|u(x, t)-\bar u(x-\alpha(t))|_{H^2}&\le
C(1+t)^{-\frac{1}{4}} |u(x,0)-\bar u(x)|_{L^1\cap H^2},\\
\dot \alpha(t) &\le C(1+t)^{-\frac{1}{2}} |u(x,0)-\bar u(x)|_{L^1\cap H^2},
\\
\alpha(t) &\le C |u(x,0)-\bar u(x)|_{L^1\cap H^2}.
\ea
Moreover, it is orbitally unstable with respect to small $H^2$ perturbations
not lying in $\cM_{cs}$, in the sense that the corresponding solution leaves
a fixed-radius neighborhood of the set of translates of $\bar u$ in
finite time.
\et


\br\label{srmks}
\textup{
Theorem \ref{t:mainstab} includes in passing the result 
that existence of unstable 
eigenvalues implies nonlinear instability, hence
completey characterizes stability/instability of waves
under the nondegeneracy condition (D1).
The rates of decay \eqref{bounds}
are exactly those of the one-dimensional case \cite{Z4}.
}
\er

\subsection{Bifurcation conditions}
We next recall the following result from \cite{TZ2,TZ3}
characterizing Hopf bifurcation of $\bar u^\eps$
in terms of conditions on the spectrum of $L(\eps)$.
Define the {\it Hopf bifurcation condition}:

\quad (D2) \quad
Outside the essential spectrum of $L(\eps)$,
for $\eps$ and $\delta>0$ sufficiently small, 
the only eigenvalues of $L(\eps)$ with real part of absolute value
less than $\delta$ are a crossing conjugate pair 
$\l_\pm (\e) := \gamma(\eps)\pm i\tau(\eps)$ of $L(\eps)$,
with $\gamma(0)=0$, $\partial_\eps \gamma(0)>0$, and $\tau(0)\ne 0$.

\begin{prop}[\cite{TZ2,TZ3}]\label{multinewPH}
Let $\bar u^\eps$, \eqref{multisysepseqn}
be a family of traveling-waves and systems 
satisfying assumptions (H0)--(H3) and (D2),
and $\eta>0$ sufficiently small.
Then, for $a\ge 0$ sufficiently small and $C > 0$ sufficiently
large, there are $C^1$ functions
$\eps(a)$, $\e(0)=0$, and $T^*(a)$, $T^*(0)=2\pi/\tau(0)$, and a $C^1$ 
family of solutions $u^{a}(x_1,t)$ of \eqref{multisysepseqn} 
with $\e=\e(a)$, time-periodic of period $T^*(a)$, such that
\be\label{abd}
   C^{-1} a \leq  \,
   \sup_{x_1 \in \R} \, e^{\eta |x_1|}
\, \big| u^{a}(x,t) - \bar u^{\e(a)}(x_1) \big| \, \leq Ca 
\qquad \hbox{\rm for all $t \geq 0$}.
\ee
Up to fixed translations in $x$, $t$, for $\eps$ sufficiently
small, these are the only nearby solutions as measured in 
norm $\|f\|_{X_1}:=\|(1+|x_1|)f(x)\|_{L^\infty(x)}$ that are
time-periodic with period
$T\in [T_0, T_1]$, for any fixed $0<T_0<T_1<+\infty$.
Indeed, they are the only 
nearby solutions of form 
$u^{a}(x,t)= {\bf u}^a(x-\sigma^a t,t)$
with $ {\bf u}^a$ periodic in its second argument.
\end{prop}

\begin{proof}
This result was established in Theorem 1.4, \cite{TZ2} with \eqref{abd}
replaced by
\be\label{xbd}
   C^{-1} a \leq  \, \sup_{x_1 \in \R} \, (1 + |x_1|) 
\, \big| u^{a}(x,t) - \bar u^{\e(a)}(x_1) \big| \, \leq Ca
\ee
and under the further assumption that there are 
no eigenvalues of $L(\eps)$ with strictly positive real part other
than possibly $\lambda_\pm(\eps)$.
As $L(\eps)$ by standard considerations \cite{He,TZ2} 
possesses at most a finite set of positive real part eigenvalues, 
an examination of the proof shows that the more general case follows
by essentially the same argument, reducing by the
Lyapunov--Schmidt reduction described in \cite{TZ2} to a finite-dimensional
equation on the direct sum of the oscillatory eigenspace
associated with $\lambda_\pm$ and the unstable eigenspace of $L$,
then appealing to standard, finite-dimensional theory to conclude
the appearance of Hopf bifurcation with bound \eqref{xbd}.
The stronger result of exponential localization, \eqref{abd},
may be obtained by combining the argument of \cite{TZ2}
with the strengthened cancellation estimates of Proposition 2.5 \cite{TZ3}.
As the distinction between \eqref{abd} and \eqref{xbd} is
not important for the present discussion, we omit the (straighforward) details.
\end{proof}

\br\label{bifchar}
\textup{
Together with Theorem \ref{t:mainstab},
Proposition \ref{multinewPH} implies that, under the Hopf bifurcation
assumption (D2) together with the further assumption that $L(\eps)$
have no strictly positive real part eigenvalues other than possibly
$\lambda_\pm$, waves $\bar u^\eps$ are linearly and nonlinearly
 {\it stable} for $\eps<0$ and {\it unstable} for $\eps>0$,
with bifurcation/exchange of stability at $\eps=0$.
}
\er

\subsection{Longitudinal vs. transverse bifurcation}
The analysis of \cite{TZ2} in fact gives slightly more
information.  Denote by
\be\label{proj}
\Pi^\e f:= \sum_{j=\pm} \phi^\eps_j(x)\langle \tilde \phi^\eps_j, f\rangle
\quad
\ee
the $L(\e)$-invariant projection
onto the oscillatory eigenspace $\Sigma^\eps:=\Span \{\phi^\eps_\pm\}$,
where $\phi^\eps_\pm$ are the eigenfunctions 
associated with $\lambda_\pm(\eps)$.  Then, we have the following result, 
proved but not explicitly stated in \cite{TZ2}.

\bpr\label{linapprox}
Under the assumptions of Proposition \ref{multinewPH},
also
\be\label{further}
\sup_{x_1}\, e^{\eta |x_1|}|u^a - \bar u-  \Pi^\eps(u^a-\bar u)|
\le Ca^2
\qquad \hbox{\rm for all $t \geq 0$}.
\ee
\epr

\begin{proof}
The weaker bound
\be\label{xfurther}
\sup_{x_1}\, (1+|x_1|)|u^a - \bar u-  \Pi^\eps(u^a-\bar u)|
\le Ca^2
\qquad \hbox{\rm for all $t \geq 0$},
\ee
is established
in the course of the Lyapunov reduction of \cite{TZ2}; see (2.17), 
Proposition 2.9, case $\omega\equiv 0$.
The stronger version \eqref{further} follows by the same argument
together with 
the strengthened cancellation estimates of Proposition 2.5 \cite{TZ3}.
\end{proof}

Bounds \eqref{abd} and \eqref{further} together
yield the standard finite-dimensional property that bifurcating solutions 
lie to quadratic order in the direction of the oscillatory eigenspace 
of $L(\eps)$.
From this, we may draw the following additional conclusions about
the structure of bifurcating waves.
By separation of variables, and $x_2$-independence of the coefficients
of $L(\eps)$, we have that the eigenfunctions $\psi$ of $L(\eps)$ decompose
into families 
\be\label{decomp}
e^{i\xi x_2}\psi(x_1),
\quad
\xi= \frac{\pi k}{M},
\ee
associated with different integers $k$, where $M$ is cross-sectional width.
Thus, there are two very different cases:
(i) ({\it longitudinal instability}) the bifurcating eigenvalues $\lambda_\pm
(\eps)$ are associated with wave-number $k=0$, or
(ii) ({\it transverse instability}) the bifurcating eigenvalues $\lambda_\pm
(\eps)$ are associated with wave-numbers $\pm k\ne 0$.

\bc
Under the assumptions of Proposition \ref{multinewPH},
$u^a$ depend nontrivially on $x_2$ if and only if the bifurcating
eigenvalues $\lambda_\pm$ are associated with transverse wave-numbers
$\pm k\ne 0$.
\ec

\begin{proof}
For $k\ne 0$, the result follows by the fact that, by
\eqref{abd} and \eqref{further}, $\Pi (u^a-\bar u)$
is the dominant part of $u^a-\bar u$, and the fact
that $\Pi f$ by inspection depends nontrivially on $x_2$
whenever $\Pi f\ne 0$.
For $k=0$, the result follows by uniqueness, and the fact
that, restricted to the one-dimensional case, the same argument
yields a bifurcating solution depending only on $x_1$.
\end{proof}

Bifurcation through longitudinal instability corresponds to
``galloping'' or ``pulsating'' instabilities described in 
detonation literature, while
symmetry-breaking bifurcation through transverse instability 
corresponds to ``cellular'' instabilities introducing 
nontrivial transverse geometry to the 
structure of the propagating wave.

\subsection{The refined stability condition and bifurcation}

Longitudinal or ``galloping'' bifurcation, 
though almost certainly occurring for
detonations (see \cite{TZ2,TZ4} and references therein), 
has up to now not been observed for shock
waves as far as we know (though we see no reason why they should
not in general be possible), nor has there been proposed any
specific mechanism by which this might occur.
The main purpose of the present paper, as we now describe, 
is to point out that for transverse or ``cellular''
bifurcations, to the contrary,
there is a simple and natural mathematical mechanism, closely related
to the inviscid stability theory for shocks in the whole space,
by which they can and likely do occur.

\subsubsection{The inviscid stability condition}
Inviscid stability analysis for shocks in the whole space
 centers about the {\it Lopatinski determinant}
\ba\label{Delta}
\Delta(\tilde \xi, \lambda)&:=
\bp
\rvec_1^- & \cdots & \rvec_{p-1}^-&
\rvec_{p+1}^+ & \cdots & \rvec_{n}^+&
\lambda [u] + i\tilde \xi [f^{2}]
\ep,\\
\ea
$\tilde \xi \in \RR^{1}$, 
$\lambda=\gamma+i\tau \in \CC$, $\tau > 0$,
a spectral determinant whose zeroes correspond to
normal modes $e^{\lambda t}e^{i\tilde \xi x_2}w(x_1)$
of the constant-coefficient linearized equations about the 
discontinuous shock solution.  Here, 
$\{ \rvec^+_{p+1}, \dots, \rvec^+_{n}\}$ and
$\{ \rvec^-_1, \dots, \rvec^-_{p-1}\}$
denote bases for the
unstable/resp. stable subspaces of 
\be\label{calA}
{\cal A}_+ (\tilde \xi,\lambda ):= 
(\lambda I+i\tilde \xi df^{2}(u_\pm))(df^1(u_\pm))^{-1}.
\ee
Weak stability $|\Delta|>0$ for $\tau>0$ is clearly necessary for linearized stability,
while strong, or uniform stability, $|\Delta|/|(\tilde \xi, \lambda)|\ge c_0>0$,
is sufficient for nonlinear stability.
Between strong instability, or failure of weak stability,
and strong stability, there lies a region of neutral stability
corresponding to the appearance of surface waves propagating
along the shock front, for which $\Delta$ is nonvanishing for
$\Re \lambda >0$ but has one or more roots
$(\tilde \xi_0, \lambda_0)$ with $\lambda_0=\tau_0$ pure imaginary.
This region of neutral inviscid stability typically occupies an open 
set in physical parameter space \cite{M1, M2, M3, BRSZ, Z1, Z2}. 
For details, see, e.g., \cite{Er1, M1, M2, M3, Me, Se1, Se2, Se3, ZS, Z1,
Z2, Z3, BRSZ}, and references therein.

It has been suggested \cite{MR1,MR2,AM} that nonlinear hyperbolic 
evolution of surface waves in the region of neutral linear stability 
might explain the onset 
of complex behavior such as Mach stem formation/kinking of the shock.
We pursue 
here a variant of this idea based instead on interaction
between neglected viscous effects and transverse spatial scales.

\subsubsection{The refined stability condition}

Viscous stability analysis for shocks in the whole space
centers about the {\it Evans function}
$$
D(\tilde \xi, \lambda),
$$
$\tilde \xi \in \RR^{1}$, 
$\lambda=\gamma+i\tau \in \CC$, $\tau > 0$,
a spectral determinant analogous to the Lopatinski determinant of the
inviscid theory, whose zeroes correspond to
normal modes $e^{\lambda t}e^{i\tilde \xi x_2}w(x_1)$,
of the linearized equations about $\bar u$ (now variable-coefficient), 
or spectra of the linearized operator about the wave.
The main result of \cite{ZS}, establishing a rigorous relation between
viscous and inviscid stability, was the asymptotic expansion
\be\label{LF}
D(\tilde \xi, \lambda)=
\gamma \Delta(\tilde \xi, \lambda) + o(|(\tilde \xi, \lambda)|)
\ee
of $D$ about the origin $(\tilde \xi, \lambda)=(0,0)$, where
$\gamma$ is a constant measuring tranversality of 
$\bar u$ as a connecting orbit of the traveling-wave ODE.
Equivalently, considering
$D(\tilde \xi,\lambda)=
D(\rho\tilde  \xi_0, \rho \lambda_0)$ as a function of polar coordinates
$(\rho,\tilde  \xi_0, \lambda_0)$, we have
\be\label{polar}
D|_{\rho=0}=0 \, \hbox{ and } \, 
(\partial/\partial \rho)|_{\rho=0} D=\gamma \Delta(\tilde \xi_0, \lambda_0).
\ee

An important consequence of \eqref{LF}
is that {\it weak inviscid stability, $|\Delta|>0$, is necessary for weak
viscous stability}, $|D|>0$ (an evident necessary condition
for linearized viscous stability).
For, \eqref{LF} implies that 
the zero set of $D$ is tangent at the origin to the
cone $\{\Delta=0\}$ (recall, \eqref{Delta}, that $\Delta$ is homogeneous,
degree one), hence enters $\{\tau>0\}$ if $\{\Delta=0\}$ does.
Moreover, in case of {\it neutral inviscid stability} 
$\Delta(\xi_0, i \tau_0)=0$, $(\xi_0, i\tau_0)\ne (0,0)$,
one may extract a further, {\it refined stability condition}
\be\label{beta}
\beta:=-D_{\rho \rho}/ D_{\rho \lambda}|_{\rho=0} \ge 0
\ee
necessary for weak viscous stability. 
For, \eqref{polar} then implies 
$D_\rho|_{\rho=0}=\gamma \Delta(\xi_0, i\tau_0)=0$, whence
Taylor expansion of $D$ yields that the zero level set of $D$
is concave or convex toward $\tau>0$ according as the sign
of $\beta$; see \cite{ZS} for details.
As discussed in \cite{ZS, Z1}, the constant $\beta$ has a heuristic
interpretation as an effective diffusion coefficient
for surface waves moving along the front.

As shown in \cite{ZS,BSZ}, the formula \eqref{beta} is well-defined 
whenever $\Delta$ is analytic at $(\tilde \xi_0, i\tau)$, in which
case $D$ considered as a function of polar coordinates
is analytic at $(0,\tilde \xi_0, i\tau_0)$, and
$i\tau_0$ is a simple root of $\Delta(\tilde \xi_0, \cdot)$.
The determinant $\Delta$ in turn is analytic at $(\tilde \xi_0,i\tau_0)$,
for all except a finite set of branch singularities
$\tau_0=\tilde \xi_0 \eta_j$.
As discussed in \cite{BSZ,Z2,Z3}, the apparently nongeneric behavior
that the family of holomorphic functions $\Delta^\eps$ associated
with shocks $(u^\eps_+,u^\eps_-)$ have roots $(\xi^\eps_0, i\tau_0(\eps))$
with $i\tau_0$ pure imaginary on an open set of $\eps$ is explained
by the fact that, on certain components of the complement on the
imaginary axis of this finite set of branch singularities,
$\Delta^\eps(\xi^\eps_0, \cdot)$ takes the imaginary axis to itself.
Thus, zeros of odd multiplicity persist on the imaginary axis, by
consideration of the topological degree of $\Delta^\eps$ as
a map from the imaginary axis to itself.

Moreover, the same topological considerations show that a simple 
imaginary root of this type can only enter or leave the imaginary axis at a 
branch singularity of $\Delta^\eps(\tilde \xi^\eps,\cdot)$ or at infinity,
which greatly aids in the computation of transition points
for inviscid stability \cite{BSZ,Z1,Z2,Z3}.
As described in \cite{Z2,Z3,Se1}, escape to infinity is always
associated with transition to strong instability.  Indeed, using
real homogeneity of $\Delta$, we may rescale by $|\lambda|$ to find
in the limit as $|\lambda|\to \infty$ that 
$
0=|\lambda_0|^{-1}\Delta(\tilde \xi_0, \lambda_0)= \Delta(\tilde \xi_0/|\lambda_0|,
\lambda_0/|\lambda_0|)\to \Delta(0,i),
$
which, by the complex homogeneity 
$\Delta(0, \lambda)\equiv \lambda\Delta(0,1)$ of the one-dimensional
Lopatinski determinant $\Delta(0, \cdot)$, yields {\it one-dimensional
instability} $\Delta(0,1)=0$.
As described in \cite{Z1}, Section 6.2,
this is associated not with surface waves, but the more dramatic
phenomenon of {\it wave-splitting}, in which the axial structure
of the front bifurcates from a single shock to a more complicated
multi-wave Riemann pattern.

\begin{exam}\label{gasrmk}
For gas dynamnics, complex symmetry,
$\bar\Delta(\tilde \xi, \lambda)= \Delta(-\tilde \xi, \bar \lambda)$,
and rotational invariance,
$\Delta(\tilde \xi, \lambda)=\Delta(-\tilde \xi, \lambda)$, imply that
\be\label{gassym}
\Delta(\tilde \xi, i\tau)= \Delta(|\tilde \xi|^2, |\tau|^2).
\ee
Explicit computation \cite{Er1,M1,Z1} yields that $\Delta(\tilde \xi_0, \cdot)$
has a pair of branch points of square-root type, located at
$|\tau_0|^2=|\tilde \xi_0|^2 (M^2-1)$, where $M$ is the downstream
Mach number $c^2/u^2$, where $c$ is sound speed and $u$ the axial
particle velocity of the shock on the downstream side, defined as
the side in the direction of particle velocity.
Transition from strong stability to neutral stability occurs
through a pair of simple imaginary zeros entering the imaginary axis at 
the branch points, and transition from neutral stability to strong instability
occurs through escape of these zeros to infinity, 
with associated one-dimensional instability/wave-splitting.
\end{exam}

\br\label{oned}
\textup{
We note in passing that one-dimensional inviscid 
stability $\Delta^\eps(0,1)\ne 0$
is equivalent to (H3) through the relation
\be\label{onedet}
\Delta^\eps(0,\lambda)=
\lambda \det(r_1^-, \dots, r_c^-, r_{c+1}^+, \dots, r_n^+, u_+-u_-).
\ee
}
\er

\subsubsection{Transverse bifurcation of flow in a duct}
We now make an elementary observation connecting
cellular bifurcation of flow in a duct to stability of shocks
in the whole space: specifically, to violation of the
refined stability condition.
Assume for the family $\bar u^\eps$ the stability conditions:

\quad (B1) \quad
For $\eps$ sufficiently small, the inviscid shock $(u^\eps_+,u^\eps_-)$
is weakly stable; more precisely, $\Delta^\eps(1, \lambda)$
has no roots $\Re \lambda \ge 0$ but a single simple pure imaginary root
$\lambda(\eps)=i\tau_*(\eps)\ne 0$ lying away 
from the singularities of $\Delta^\eps$.

\medskip
\quad (B2) \quad 
The refined stability coefficient $\beta(\eps)$
defined in \eqref{beta} satisfies $\Re \beta(0)=0$, $\partial_\eps \Re \beta(0)<0$.

\medskip

\bl[\cite{ZS,Z1}]\label{branch}
Assuming (H0)--(H2) and (B1), for $\eps$, $\tilde \xi$ sufficiently
small, there exist a smooth family of roots 
$(\tilde \xi, \lambda_*^\eps(\tilde \xi))$ of $D(\tilde \xi, \lambda)$
with
\be\label{expansion}
\lambda_*^\eps(\tilde \xi)= i\tilde \xi \tau_*(\eps) 
- \tilde \xi^2 \beta(\eps) +\delta (\eps)\tilde \xi^3 + r(\eps, \tilde \xi)
\tilde \xi^4,
\qquad
r\in C^1(\eps,\tilde \xi). 
\ee
Moreover, these are the unique roots of $D$ satisfying
$\Re \lambda \ge -\tilde |\xi|/C$ for some $C>0$
and $\rho=|(\tilde \xi, \lambda)|$ sufficiently small.
\el

\begin{proof}
This follows by the Implicit Function Theorem applied
to the function $\check D(\rho, \lambda_0):=\rho^{-1}
D(\rho \tilde \xi_0, \rho \lambda_0)$, 
about the values $(\rho, \tilde \lambda_0)=(0,i\tilde \xi_0 \tau_*(\eps)$,
where $\tilde \xi_0$ is without loss of generality held fixed,
using the facts that $D$ expressed in
polar coordinates $(\rho, \tilde \xi_0, \lambda_0)$
satisfies $D|_{\rho=0}\equiv 0$ and
$\partial_\rho D|_{\rho=0}\equiv \Delta$,
so that $D_\rho$, $D_{\lambda \lambda}$, and $D_\lambda$
all vanish at $(0, \tilde \xi_0, \tilde \xi_0 i\tau_*(\eps))$.
For details, see the proof of Theorem 3.7, \cite{Z1}.
\end{proof}

\bc\label{imag}
Assuming (H0)--(H2) and (B1)-(B2), for $\eps$, $\tilde \xi$ sufficiently
small, there is a unique $C^1$ function $\cE(\tilde \xi)\ne 0$,
$\cE(0)=O$,
such that $\Re \lambda_*^\eps(\tilde \xi)=0$ for $\eps=\cE(
\tilde \xi)$. In the generic case $\Re \delta(0)\ne 0$, 
moreover, 
\be\label{param}
\cE(\tilde \xi) \sim (\Re \delta(0)/\partial_\eps \beta(0)) \tilde \xi.
\ee
\ec

\begin{proof}
As a consequence of \eqref{expansion}, we have for some smooth $G$
\be\label{realexp}
\Re \Big(\frac{\lambda_*^\eps(\tilde \xi)}{\tilde \xi^2}\Big)= 
- \Re \beta(\eps) + \tilde \xi G(\eps, \tilde \xi),
\ee
$G:=(\delta(\eps)+r(\eps, \tilde \xi))$,
whence the equation
$0=\Re \Big(\frac{\lambda_*^\eps(\tilde \xi)}{\tilde \xi^2}\Big)= 
- \Re \beta(\eps) + \tilde \xi G(\tilde \xi, \eps)$ has a unique
root $\eps=\cE(\tilde \xi)$ by assumption (B2) and standard
scalar bifurcation theory.
From $G=\delta(\eps)\tilde \xi + O(|\tilde \xi|^2)$, we find,
in the generic case $\Re \delta(0)\ne 0$, 
that $\partial_{\tilde \xi}(\tilde \xi G)|_{\tilde \xi, \eps=0,0}=
\Re \delta(0)\ne 0$, yielding \eqref{param}.
\end{proof}

\br\label{realapp}
\textup{
As described further in \cite{BSZ}, the above results on the refined stability
condition apply also in the case of ``real'' or partial viscosity,
in particular to the physical Navier--Stokes equations of compressible 
gas dynamics and MHD.
}
\er

To (B1) and (B2), adjoin now the additional assumptions:

\quad (B3) \quad $\delta(0)\ne 0$.
\medbreak

\quad (B4) \quad 
At $\eps=0$, the Evans function $D(\tilde \xi, \lambda)$ has
no roots $\tilde \xi\in \RR$, $\Re \lambda \ge 0$ outside 
a sufficiently small ball about the origin.
\medbreak

Then, we have the following main result. 

\bt\label{refbif}
Assuming (H0)--(H2) and (B1)-(B4), for
$\eps_{max}>0$ sufficiently small and each
cross-sectional width $M$ sufficiently large,
there is a finite sequence
$0<\eps_1(M)<\dots < \eps_k(M)<\dots \le \eps_{max}$,
with $\eps_k(M) \sim (\Re \delta(0)/\partial_\eps \beta(0)) 
\frac{\pi k}{M}$,
such that, as $\eps$ crosses successive $\eps_k$ from the left, there occur
a series of transverse 
(i.e., ``cellular'') Hopf bifurcations of $\bar u^\eps$ 
associated with wave-numbers $\pm k$,
with successively smaller periods $T_k(\eps) \sim \tau_*(0) \frac{2M}{k}$.
\et

\begin{proof}
By (B4), for $|\eps|\le \eps_{max}$ sufficiently small,
we have by continuity that there exist no roots of $D(\tilde \xi, \lambda)$
for $\Re \lambda\ge -1/C$, $C>0$, outside a small ball about the origin.
By Lemma \ref{branch}, within this small ball, there are no roots
other than possibly $(\tilde \xi, \lambda^\eps(\tilde \xi))$
with $\Re \lambda\ge 0$: in particular, 
{\it no nonzero purely imaginary spectra are possible 
other than at values $\lambda^\eps(\tilde \xi)$} for 
operator $L(\eps)$ acting on functions on the whole space.

Considering $L$ instead as an operator acting on functions 
on the channel $\CalC:= \{x:\, (x_1, x_2)\in \RR^1\times [-M,M]\}$,
we find by discrete Fourier transform/separation of variables
that its spectra are exactly the zeros of $D(\xi_k,\lambda)$,
as $\xi_k= \frac{\pi k}{L}$ runs through all integer
wave-numbers $k$; see \eqref{decomp}.
Applying Corollary \ref{imag}, and using (B3), we find, therefore,
that pure imaginary eigenvalues of $L(\eps)$ with $|\eps|\le \eps_{max}$
sufficiently small occur precisely at values $\eps=\eps_k$, and
consist of crossing conjugate pairs $\lambda_\pm^k(\eps)$
associated with wave-numbers $\pm k$, 
satisfying Hopf bifurcation condition (D2) with
$$
\Im \lambda_\pm^k(\eps)\sim   \tau_*(0)\Big)\frac{\pi k}{L}.
$$
Applying Proposition \ref{multinewPH}, we obtain the result.
\end{proof}

\br\label{doubling}
\textup{
As evidenced by decreasing periods $T_k$, this phenomenon
of increasing-complexity solutions
is completely different from the more familiar one of period-doubling.
}
\er

\br\label{B1alt}
\textup{
Lemma \ref{branch} and Corollary \ref{imag} are readily generalized
to the case with (B1) replaced by  
(B1') For $\eps$ sufficiently small, the inviscid shock $(u^\eps_+,u^\eps_-)$
is weakly stable, with all pure imaginary roots simple and
lying away from the singularities of $\Delta^\eps$.
In this case we obtain
a famiy of roots/crossings of the imaginary axis, one for each 
imaginary root of $\Delta^\eps$.
In particular, for gas dynamics, due to rotational invariance
(see example \ref{gasrmk}), we obtain families of four crossing
eigenvalues $\lambda_\pm(\eps_k)$, with each of $\lambda_+$ and
$\lambda_-$ occurring at both wave-numbers $k$ and $-k$.
This is not a standard Hopf bifurcation, but a more complicated
version with $O(2)$ symmetry, and so we cannot apply directly
Theorem \ref{refbif}.
}
\er

\subsection{Discussion and open problems}
We have presented in a simple setting a rigorous mathematical
demonstration of a mechanism by which destabilization
of hyperbolic surface waves arising in the inviscid shock stability
problem in the whole space
can, at appropriate transverse length scales, lead to Hopf bifurcation
of a viscous shock in a finite-cross-section duct: specifically,
destabilization of the effective transverse viscosity $\beta$
investigated in \cite{ZS,BSZ}, or violation of the {\it refined
stability condition}.
This appears to be a fundamentally different mechansim than the 
hyperbolic ones proposed by Majda et al \cite{MR1,MR2,AM} via weakly noninear
geometric optics.

We point out that as shock parameters cross the inviscid strong
instability boundary, gas-dynamical shocks undergo one-dimensional
instability, or {\it wave-splitting}, 
a more dramatic change in front topology than the cellular
instabilities we seek to investigate.
Thus, cellular instability must occur {\it before} the strong
instability boundary is reached.
Experimental observations, though not conclusive, indicate
that nonetheless it occurs near the strong instability boundary
\cite{BE}, suggesting that it lies in the region of neutral
inviscid stability as we have conjectured.

More,
if the transition to cellular instability occurs at low frequencies,
it {\it must} occur by the scenario described, or else the viscous
shock would remain stable up to the point of wave-splitting.
If, on the other hand, it occurs at high frequencies, then as
pointed out in \cite{BSZ}, then it necessarily involves Hopf
bifurcation, by one-dimensional inviscid stability, (H2)
(satisfied for typical equations of state).
Thus, it would appear quite promising to search for Hopf bifurcations
in the region of neutral inviscid stability, whether of the 
``low-frequency'' type studied here or a ``high-frequency'' type
involving unknown mechanisms.
This would be a very interesting direction for numerical investigations,
for example by the numerical Evans function techniques of
\cite{Br1,Br2,BrZ,BDG,HuZ,BHRZ}.

Another interesting direction would be investigation
of the stability coefficient $\beta$, both numerically and analytically.
As pointed out in \cite{BSZ}, this is numerically quite well-conditioned.
One might also consider attempting to carry out an asymptotic analysis
near the endpoints of the region of neutral stability, at which the
imaginary root $\tau_*$ approaches either a branch point of $\Delta$
or else infinity.

At a technical level, an interesting open problem
is to carry out a bifurcation analysis in the rotationally
symmetric case, for example, for gas dynamics,
in which the bifurcation associated with crossing $\lambda_\pm$
no longer a standard Hopf bifurcation
but a more complicated type involving $O(2)$ symmetry.
For a description of Hopf bifurcation with $O(2)$ symmetry, see, for
example, \cite{W}.
For a circular cross-section, there is besides
axial translation an additional continuous
group invariance of rotation in the transverse direction,
leading to the possibility of ``spinning'' instabilities.
These degenerate cases require further analysis at the level
of the finite-dimensional reduced equations; 
however, the initial reduction to finite dimensions, as carried
out in \cite{TZ2}, is essentially the same.
Other open problems are to
extend to detonations, as done for the one-dimensional case in \cite{TZ4} 
and to treat also real, or partial viscosities.
As discussed in \cite{TZ2,TZ3}, the latter problem involves interesting
issues involving Lagrangian vs. Eulerian formulations.


\section{Conditional stability analysis}\label{s:cond}
Nonlinear stability follows quite similarly
as in the one-dimensional case \cite{Z4}, decomposing
behavior into a one-dimensional (averaged in $x_2$) flow
driving time-exponentially damped transverse modes.\footnote{
Indeed, though we do not do it here, this prescription could
be followed quite literally at the nonlinear level.}

Define the perturbation variable
\be\label{pert2}
v(x,t):=u(x+\alpha(t),t)-\bar u(x)
\ee
for $u$ a solution of \eqref{multiNS}--\eqref{BC}, 
where $\alpha$ is to be specified later.
Subtracting the equations for $u(x+\alpha(t), t)$ and $\bar u(x)$,
we obtain the nonlinear perturbation equation
\be\label{nlpert}
v_t-Lv= \sum_{j=1}^2 N_j(v)_{x_j} +\d_t \alpha (\bar u_{x_1}+ \partial_{x_1} v),
\ee
where 
\be\label{L}
L:=\Delta_x - \sum_{j=1}^2  \partial_{x_j}A^j(x),
\quad A^j:=df_j(\bar u)
\ee
denotes the linearized operator about $\bar u$ and
$
N_j(v):=-(f^j(\bar u+ v)-f^j(\bar u)- df^j(\bar u)v),
$
where, so long as $|v|_{H^1}$ (hence $|v|_{L^\infty}$ and $|u|_{L^\infty}$) 
remains bounded, 
\ba\label{Nbds}
N^j(v)&=O(|v|^2),
\qquad
\partial_x N^j(v)=O(|v||\partial_x v|),
\qquad
\partial_x^2 N^j(v)= O(|\partial_x|^2+ |v||\partial_x^2v|).\\
\ea

\subsection{Projector bounds}\label{projbds}
Let $\Pi_u$ denote the eigenprojection of $L$ onto its
unstable subspace $\Sigma_u$, and $\Pi_{cs}={\rm Id}- \Pi_u$
the eigenprojection onto its center stable subspace $\Sigma_{cs}$. 

\bl\label{projlem}
Assuming (H0)--(H1), there is $\tilde \Pi_j$ defined in \eqref{tildePidef}
such that
\be\label{comm}
\Pi_j \partial_x= \partial_x \tilde \Pi_j
\ee
for $j=u,\, cs$ and, for all $1\le p\le \infty$, $0\le r\le 4$, 
\ba\label{pbd}
|\Pi_{u}|_{L^{p}\to W^{r,p}}, |\tilde \Pi_{u}|_{L^{p}\to W^{r,p}}&\le C,\\
|\tilde \Pi_{cs}|_{W{r,p}\to W^{r,p} }, 
\; |\tilde \Pi_{cs}|_{W{r,p}\to W^{r,p}}&\le C.
\ea
\el

\begin{proof}
Recalling that $L$ has at most finitely many unstable eigenvalues,
we find that $\Pi_u$ may be expressed as
$$
\Pi_u f= \sum_{j=1}^p \phi_j(x) \langle \tilde \phi_j, f\rangle,
$$
where $\phi_j$, $j=1, \dots p$ are generalized right eigenfunctions of
$L$ associated with unstable eigenvalues $\lambda_j$, 
satisfying the generalized eigenvalue equation $(L-\lambda_j)^{r_j}\phi_j=0$,
$r_j\ge 1$, and $\tilde \phi_j$
are generalized left eigenfunctions.
Noting that $L$ is divergence form, and that $\lambda_j\ne 0$,
we may integrate $(L-\lambda_j)^{r_j}\phi_j=0$ over $\R$ to
obtain $\lambda_j^{r_j}\int \phi_j dx=0$ and thus $\int\phi_j dx=0$.
Noting that $\phi_j$, $\tilde \phi_j$ and derivatives decay exponentially
in $x_1$ by separation of variables and standard one-dimensional
theory \cite{He,ZH,MaZ1}, we find that
$
\phi_j= \partial_x \Phi_j
$
with $\Phi_j$ and derivatives exponentially decaying in $x_1$, hence
\be\label{tildePidef}
\tilde \Pi_u f=\sum_j \Phi_j \langle \partial_x \tilde \phi, f\rangle.
\ee
Estimating 
$
|\partial_x^j\Pi_u f|_{L^p}=|\sum_j \partial_x^j\phi_j \langle \tilde \phi_j f
\rangle|_{L^p}\le
\sum_j |\partial_x^j\phi_j|_{L^p} |\tilde \phi_j|_{L^q} |f|_{L^p}
\le C|f|_{L^p}
$
for $1/p+1/q=1$
and similarly for $\partial_x^r \tilde \Pi_u f$, we obtain the claimed
bounds on $\Pi_u$ and $\tilde \Pi_u$, from which the bounds on
$\Pi_{cs}={\rm Id}-\Pi_u$ and
$\tilde \Pi_{cs}={\rm Id}-\tilde \Pi_u$ follow immediately.
\end{proof}

\subsection{Linear estimates}

Let $G_{cs}(x,t;y):=\Pi_{cs}e^{Lt}\delta_y(x)$ denote
the Green kernel of the linearized solution operator on
the center stable subspace $\Sigma_{cs}$.
Then, we have the following detailed pointwise bounds
established in \cite{TZ2,MaZ1}.

\begin{proposition}[\cite{TZ2,MaZ1}]\label{greenbounds}
Assuming (H0)--(H2), (D1)--D(3), 
the center stable Green function may be decomposed as 
$G_{cs}=E+\tilde G$, where 
\begin{equation}\label{E}
E(x,t;y)= \d_{x_1} \bar u(x_1) e(y_1,t),
\end{equation}
\begin{equation}\label{e}
  e(y_1,t)=\sum_{a_k^{-}>0}
  \left(\textrm{errfn }\left(\frac{y_1+a_k^{-}t}{\sqrt{4(t+1)}}\right)
  -\textrm{errfn }\left(\frac{y_1-a_k^{-}t}{\sqrt{4(t+1)}}\right)\right)
  l_{k}^{-}(y_1)
\end{equation}
for $y_1\le 0$ and symmetrically for $y_1\ge 0$, 
$l_k^-\in \R^n$ constant, and
$a_j^\pm$ are the eigenvalues of $df(u_\pm)$, and

\be \label{tGbounds}
| \int_{\CalC} \d_x^s \tG(\cdot,t;y)f(y)dy|_{L^p}
\le C (1+t^{-\frac{s}{2}})t^{-\frac{1}{2}(\frac{1}{q}-\frac{1}{p})} |f|_{L^q},
\ee
\be\label{tGybounds}
|
\int_{\CalC} \d_x^s\tG_y(\cdot,t;y)f(y)dy|_{L^p}
\le C(1+t^{-\frac{s}{2}}) t^{-\frac{1}{2}(\frac{1}{q}-\frac{1}{p})-\frac{1}{2}} |f|_{L^q},
\ee
for all $t\ge 0$, $0\le s\le 2$, some $C>0$, for any
$1\le q\le p$ 
and $f\in L^q\cap L^p$.
\end{proposition}

\begin{proof}
As observed in \cite{TZ2},
it is equivalent to establish decomposition 
\be\label{fulldecomp}
G=G_u + E+\tilde G
\ee
for the full Green function $G(x,t;y):=e^{Lt}\delta_y(x)$,
where 
$$
G_u(x,t;y):=\Pi_u e^{Lt}\delta_y(x)
=
e^{\gamma t}\sum_{j=1}^p\phi_j(x)\tilde \phi_j(y)^t
$$
for some constant matrix $M\in \C^{p\times p}$
denotes the Green kernel of the linearized solution operator
on $\Sigma_u$, $\phi_j$ and $\tilde\phi_j$ right and left
generalized eigenfunctions associated with unstable eigenvalues
$\lambda_j$, $j=1,\dots,p$.

Using separation of variables, moreover, we may decompose
$G=\sum_k G^k$, where $G^k$ is the Green function acting
on Fourier modes of wave number $k$, i.e.,
$G^k=\CalF^{-1} G \delta(\tilde \xi- \pi k/M) \CalF$, 
where $\CalF$ denotes Fourier transform in $x_2$, and
$\tilde \xi$ Fourier frequency.
This reduces the problem to that of deriving the asserted bounds
separately on the one-dimensional Green function $G^0$
and on the complement $\sum_{k\ne 0}G^k$, where the difficulty,
due to lack of spectral gap,
is concentrated in the estimation of the one-dimensional part $G^0$.

The bounds on the one-dimensional Green function $G^0$ have already been
established in \cite{Z4}, Proposition 4.2, by essentially the same stationary 
phase estimates used in \cite{MaZ3} in the stable case $\Pi_u=0$; 
see \cite{TZ2,Z4} for further discussion.
The bounds on the complement $\sum_{k\ne 0}G^k$ follow by 
the straightforward semigroup estimate 
$|e^{\tilde Lt}f|_{L^p}\le Ce^{-\eta t} |f|_{L^p}$, $\eta>0$, where
$\tilde L$ denotes the projection of $L$ onto the intersection
of its center stable subspace and the subspace of functions with
transverse Fourier wave numbers $\ne 0$, which evidently has a nonzero
spectral gap $\sigma (\tilde L)\le -2\eta<0$ for some $\eta>0$;
see \cite{TZ2} for related computations.
\end{proof}

\bc[\cite{Z4}]\label{ebds}
The kernel ${e}$ satisfies for all $t>0$
$$
|{e}_y (\cdot, t)|_{L^p},  |{e}_t(\cdot, t)|_{L^p} 
\le C t^{-\frac{1}{2}(1-1/p)},
\label{36}
$$
$$
|{e}_{ty}(\cdot, t)|_{L^p} 
\le C t^{-\frac{1}{2}(1-1/p)-1/2}.
\label{37}
$$
\ec

\begin{proof}
Direct computation using definition \eqref{e}.
\end{proof}


\subsection{Reduced equations}
Recalling that $\d_{x_1}\bar u$ is a stationary
solution of the linearized equations $u_t=Lu$,
so that $L\d_{x_1}\bar u =0$, or
$$
\int_{\CalC} G(x,t;y)\bar u_{x_1}({y_1})dy=e^{Lt}\bar u_{x_1}({x_1})
=\d_{x_1}\bar u(x_1),
$$
we have, applying Duhamel's principle to \eqref{nlpert},
$$
\begin{array}{l}
  \displaystyle{
  v(x,t)=\int_{\CalC} G(x,t;y)v_0(y)\,dy } \\
  \displaystyle{\qquad
  -\int^t_0 \int_{\CalC} G_y(x,t-s;y)
  (N(v)+\dot \alpha v ) (y,s)\,dy\,ds + \alpha (t)\d_{x_1} \bar u({x_1}).}
\end{array}
$$
Defining 
\begin{equation}
 \begin{array}{l}
  \displaystyle{
  \alpha (t)=-\int_{\CalC}e(y,t) v_0(y)\,dy }\\
  \displaystyle{\qquad
  +\int^t_0\int_{\CalC} e_{y}(y,t-s)(N(v)+
  \dot \alpha\, v)(y,s) dy ds, }
  \end{array}
 \label{alpha}
\end{equation}
following \cite{ZH,Z4,MaZ2,MaZ3}, 
where $e$ is defined as in \eqref{e}, 
and recalling the decomposition $G=E+ G_u+ \tilde G$ of \eqref{fulldecomp},
we obtain the {\it reduced equations}
\begin{equation}
\begin{array}{l}
 \displaystyle{
  v(x,t)=\int_{\CalC} (G_u+\tilde G)(x,t;y)v_0(y)\,dy }\\
 \displaystyle{\qquad
  -\int^t_0\int_{\CalC}(G_u+\tilde G)_y(x,t-s;y)(N(v)+
  \dot \alpha v)(y,s) dy \, ds, }
\end{array}
\label{v}
\end{equation}
and, differentiating (\ref{alpha}) with respect to $t$,
and observing that 
$e_y (y_1,s)\to 0$ as $s \to 0$, as the difference of 
approaching heat kernels,
\begin{equation}
 \begin{array}{l}
 \displaystyle{
  \dot \alpha (t)=-\int_{\CalC}e_t(y,t) v_0(y)\,dy }\\
 \displaystyle{\qquad
  +\int^t_0\int_{\CalC} e_{yt}(y,t-s)(N(v)+
  \dot \alpha v)(y,s)\,dy\,ds. }
 \end{array}
\label{alphadot}
\end{equation}
\medskip

\subsection{Nonlinear damping estimate}

\begin{proposition}[\cite{MaZ3}]\label{damping}
Assuming (H0)-(H3), let $v_0\in H^{2}$, 
and suppose that for $0\le t\le T$, the $H^{2}$ norm of $v$
remains bounded by a sufficiently small constant, for $v$ as in
\eqref{pert2} and $u$ a solution of \eqref{multiNS}--\eqref{BC}.
Then, for some constants $\theta_{1,2}>0$, for all $0\leq t\leq T$,
\begin{equation}\label{Ebounds}
\|v(t)\|_{H^2}^2 \leq C e^{-\theta_1 t} \|v(0)\|^2_{H^2} 
+ C \int_0^t e^{-\theta_2(t-s)} (|v|_{L^2}^2 + |\dot \alpha|^2) (s)\,ds.
\end{equation}
\end{proposition}

\begin{proof}
Energy estimates identical with those of the one-dimensional proof in \cite{Z4},
using the fact that boundary terms in $x_2$ are identically zero
due to periodic boundary conditions.
\end{proof}

\subsection{Proof of nonlinear stability}
Decompose the nonlinear perturbation $v$ as
\be\label{vdecomp}
v(x,t)=w(x,t)+z(x,t),
\ee
where
\be\label{wzdef}
w:=\Pi_{cs}v, \quad z:=\Pi_u v.
\ee
Applying $\Pi_{cs}$ to \eqref{v} and recalling commutator
relation \eqref{comm}, we obtain an equation
\ba \label{w}
  w(x,t)&=\int_{\CalC} \tilde G (x,t;y)w_0(y)\,dy \\
  &\quad -\int^t_0\int_{\CalC} \tilde G_y (x,t-s;y)
\tilde \Pi_{cs}(N(v)+
  \dot \alpha v)(y,s) dy \, ds
\ea
for the flow along the center stable manifold, parametrized by
$w\in \Sigma_{cs}$.

\bl\label{quadlem} Assuming (H0)--(H1), for $v$ lying initially
on the center stable manifold $\cM_{cs}$,
\be\label{vwbd}
|z|_{W^{r,p}}\le C|w|_{H^2}^2
\ee 
for some $C>0$, for all $1\le p\le \infty$ and $0\le r\le 4$,
so long as $|w|_{H^2}$ remains sufficiently small.
\el

\begin{proof}
By tangency of the center stable manifold to $\Sigma_{cs}$, we
have immediately $|z|_{H^2}\le C|w|_{H^2}^2$, whence
\eqref{vwbd} follows by equivalence of norms for finite-dimensional
vector spaces, applied to the $p$-dimensional subspace $\Sigma_u$.
(Alternatively, we may see this by direct computation using
the explicit description of $\Pi_u v$ afforded by Lemma \ref{projlem}.)
\end{proof}

\begin{proof}[Proof of Theorem \ref{t:mainstab}]

Recalling by Theorem \ref{t:maincs}
that solutions remaining for all time in a sufficiently
small radius neighborhood $\cN$ of the set of translates of $\bar u$
lie in the center stable manifold $\cM_{cs}$, we obtain trivially
that solutions not originating in $\cM_{cs}$ must exit $\cN$ in finite time,
verifying the final assertion of orbital instability with respect
to perturbations not in $\cM_{cs}$.

Consider now a solution $v \in \cM_{cs}$, or, equivalently,
a solution $w\in \Sigma_{cs}$ of \eqref{w} with 
$z=\Phi_{cs}(w)\in \Sigma_u$.
Define
\begin{equation}
\label{zeta2}
 \zeta(t):= \sup_{0\le s \le t}
 \Big( |w|_{H^2}(1+s)^{\frac{1}{4}} + 
(|w|_{L^\infty}+ |\dot \alpha (s)|) (1+s)^{\frac{1}{2}} \Big).
\end{equation}
We shall establish:

{\it Claim.} For all $t\ge 0$ for which a solution exists with
$\zeta$ uniformly bounded by some fixed, sufficiently small constant,
there holds
\begin{equation}
\label{claim}
\zeta(t) \leq C_2(E_0 + \zeta(t)^2)
\quad \hbox{\rm for} \quad
E_0:=|v_0|_{L^1\cap H^2}.
\end{equation}
\medskip

{}From this result, provided $E_0 < 1/4C_2^2$, 
we have that $\zeta(t)\le 2C_2E_0$ implies
$\zeta(t)< 2C_2E_0$, and so we may conclude 
by continuous induction that
 \begin{equation}
 \label{bd}
  \zeta(t) < 2C_2E_0
 \end{equation}
for all $t\geq 0$, 
from which we readily obtain the stated bounds.
(By standard short-time $H^s$ existence theory, 
$v\in H^2$ exists and $\zeta$ remains
continuous so long as $\zeta$ remains bounded by some uniform constant,
hence \eqref{bd} is an open condition.)
\medskip

{\it Proof of Claim.}
By \eqref{pbd}, $|w_0|_{L^1\cap H^2}=|\Pi_{cs}v_0|_{L^1\cap H^2}\le CE_0$.
Likewise, by Lemma \ref{quadlem},
\eqref{zeta2}, \eqref{Nbds}, and Lemma \ref{projlem}, for $0\le s\le t$,
\ba\label{Nlast}
|\tilde \Pi_{cs}(N(v)+ \dot \alpha v)(y,s)|_{L^2}&\le C\zeta(t)^2 
(1+s)^{-\frac{3}{4}}.\\
\ea

Combining the latter bounds with representations \eqref{w} and \eqref{alphadot}
and applying Proposition \ref{greenbounds}, we obtain
 \ba\label{claimw}
  |w(x,t)|_{L^p} &\le
  \Big|\int_{\CalC} \tilde G(x,t;y) w_0(y)\,dy\Big|_{L^p}
   \\
 &\qquad +
\Big|\int^t_0
  \int_{\CalC} \tilde G_y(x,t-s;y)  \tilde \Pi_{cs}(N(v)+
  \dot \alpha v)(y,s)  dy \, ds\Big|_{L^p} \\
  & \le
  E_0 (1+t)^{-\frac{1}{2}(1-\frac{1}{p})}
+ C\zeta(t)^2 \int^t_0
(t-s)^{-\frac{3}{4}+\frac{1}{2p}}
(1+s)^{-\frac{3}{4}} dy \, ds \\
&\le
C(E_0+\zeta(t)^2)(1+t)^{-\frac{1}{2}(1-\frac{1}{p})}
\ea
and, similarly, using H\"older's inequality and
applying Corollary \ref{ebds},
\ba\label{claimalpha}
 |\dot \alpha(t)| &\le \int_{\CalC}|e_t(y,t)|
  |v_0(y)|\,dy \\
  &\qquad +\int^t_0\int_{\CalC} |e_{yt}(y,t-s)||
(N(v)+ \dot \alpha v)(y,s)|\,dy\,ds\\
&\le |e_t|_{L^\infty} |v_0|_{L^1}
+ C\zeta(t)^2 \int^t_0
|e_{yt}|_{L^2}(t-s) |(N(v)+ \dot \alpha v)|_{L^2}(s) ds\\
&\le E_0 (1+t)^{-\frac{1}{2}}
+ C\zeta(t)^2 \int^t_0
(t-s)^{-\frac{3}{4}}(1+s)^{-\frac{3}{4}} ds\\
&\le C(E_0+\zeta(t)^2)(1+t)^{-\frac{1}{2}}.\\
\ea

By Lemma \ref{quadlem}, 
\be\label{zh2}
|z|_{H^2}(t)\le C|w|_{H^2}^2(t)\le C\zeta(t)^2.
\ee
In particular, 
$
|z|_{L^2}(t)\le C\zeta(t)^2(1+t)^{-\frac{1}{2}}.
$
Applying Proposition \ref{damping} and using \eqref{claimw} and
\eqref{claimalpha}, we thus obtain
\be\label{claimwH2}
|w|_{H^2}(t)\le C(E_0+\zeta(t)^2)(1+t)^{-\frac{1}{4}}.
\ee
Combining \eqref{claimw}, \eqref{claimalpha}, and \eqref{claimwH2},
we obtain \eqref{claim} as claimed.
As discussed earlier,
from \eqref{claim}, we obtain by continuous induction \eqref{bd}, or
$
 \zeta\le 2C_2|v_0|_{L^1\cap H^2}, 
$
whereupon the claimed bounds on $|v|_{L^p}$ and $|v|_{H^2}$ follow by
\eqref{claimw} and \eqref{claimwH2}, and on $|\dot \alpha|$ 
by \eqref{claimalpha}.
Finally, a computation parallel to \eqref{claimalpha} 
(see, e.g., \cite{MaZ3,Z2}) 
yields $| \alpha(t)| \le C(E_0+\zeta(t)^2)$, from which
we obtain the last remaining bound on $|\alpha(t)|$.
\end{proof}




%


\end{document}